# FOCUSING ON BEST PROXIMITY POINTS OF GENERALIZED VERSIONS OF CYCLIC IMPULSIVE SELF-MAPPINGS

M. De la Sen

Institute of Research and Development of Processes. University of Basque Country

Campus of Leioa (Bizkaia) - Aptdo. 644- Bilbao, 48080- Bilbao. SPAIN

email: *manuel.delasen@ehu.es*

**Abstract**. This paper is devoted to the study of convergence properties of distances between points and the existence and uniqueness of best proximity and fixed points of the so-called semi-cyclic impulsive self-mappings on the union of a number of nonempty subsets in metric spaces. The concept of semi-cyclic self- mappings generalizes the well-known of cyclic ones in the sense that the iterated sequences built through such mappings are allowed to have images located in the same subset as its pre-image. The self-mappings under study may be impulsive since eventually being composite mappings involving two self-mappings, one of them being eventually discontinuous so that the formalism can potentially be applied to the study of stability of a class of impulsive differential equations and their discrete counterparts. Some application examples are also given.

## 1. Introduction

Fixed Point Theory has an increasing interest in research in the last years especially because of its high richness in bringing together several fields of Mathematics including classical and functional Analysis, topology and geometry, [1-8]. There are many fields for the potential application of this rich theory in Physics, Chemistry and Engineering, as for instance, because of its usefulness for the study of existence, uniqueness and stability of the equilibrium points and for the study of the convergence of state- solution trajectories of differential/difference equations and continuous, discrete, hybrid and fuzzy dynamic systems as well as the study of the convergence of iterates associated to the solutions. A basic key point in this context is that fixed points are equilibrium points of solutions of most of many of the above problems. Fixed Point Theory has also been investigated in the context of the so-called cyclic self-mappings, [8-13]. One of the relevant problems under study in Fixed Point Theory is that associated with $p-$ cyclic mappings which are defined on the union of a number of nonempty subsets $A_i \subset X$; $\forall i \in \bar{p} = \{1, 2, \ldots, p\}$ of metric $(X, d)$ or Banach spaces $(X, \|\ \|)$. There is an abundant background literature concerning non-expansive, non-spreading and contractive $p-$ cyclic self-mappings $T : \bigcup_{i \in \bar{p}} A_i \to \bigcup_{i \in \bar{p}} A_i$ A key point in the study of contractive cyclic self-mappings is that if the subsets $A_i$ for $i \in \bar{p}$ are disjoint then the convergence of the sequence of iterates $x_{n+1} = Tx_n$; $\forall n \in \mathbf{Z}_{0+}$ $(\mathbf{Z}_{0+} = \mathbf{Z}_+ \cup \{0\})$, $x_0 \cup \bigcup_{i \in \bar{p}} A_i$ is only possible to best proximity points. The existence of such fixed points, its uniqueness and associated properties is studied rigorously in [11-13] in the framework of uniformly convex metric spaces. This paper is focused on the study of the properties of iterated sequences $x_{n+1} = Tx_n$; $\forall n \in \mathbf{Z}_{0+}$ for any given $x_0 \in \bigcup_{i \in \bar{p}} A_i$ generated from nonexpansive and, in particular, contractive $p-$ semi-cyclic impulsive self-mappings $T(\equiv T^+ \circ T^-) : \bigcup_{i \in \bar{p}} A_i \to \bigcup_{i \in \bar{p}} A_i$, where $A_i \subset X$; $\forall i \in \bar{p}$ and $(X, d)$ is a metric space, and

(1) $T^-: \bigcup_{i\in \bar{p}} A_i \to \bigcup_{i\in \bar{p}} A_i$ is a nonexpansive, or contractive, $p$-semi-cyclic self-mapping in the sense that $T(A_i) \subseteq A_i \cup A_{i+1}$ ; $\forall i \in \bar{p}$ (being, in particular) a $p$-cyclic self-mapping if $T(A_i) \subseteq A_{i+1}$ ; $\forall i \in \bar{p}$ ); and

(2) $T^+: \bigcup_{i\in \bar{p}} A_i \to \bigcup_{i\in \bar{p}} A_i$ satisfies a condition of the type $d(T^+(T^- x), T^+(T^- y)) \leq m(T^- x, T^- y) d(T^- x, T^- y)$ for $x, y \in \bigcup_{i\in \bar{p}} A_i$ and some given bounded function $m: (\bigcup_{i\in \bar{p}} A_i) \times (\bigcup_{i\in \bar{p}} A_i) \to \mathbf{R}_{0+}$.

The properties of boundedness and convergence of distances are studied in metric spaces while those of the iterated sequences $x_{n+1} = T x_n$; $\forall n \in \mathbf{Z}_{0+}$, are studied in uniformly convex Banach spaces. On the other hand, the boundedness of the sequences of distances between consecutive iterates is guaranteed for nonexpansive $p$-semi-cyclic self-mappings while its convergence is proved for asymptotically contractive $p$-semi-cyclic self-mappings. In this case, it is proven that a limit set of such sequences exists which contains best proximity points if the asymptotically contractive $p$-semi-cyclic self-mapping is asymptotically $p$-cyclic, $(X, d)$ is a complete metric space which is also a uniformly convex Banach space $(X, \|\ \|)$ and the subsets $A_i \subset X$ ; $\forall i \in \bar{p}$ are nonempty, closed and convex.

## 2. Nonexpansive and contractive semi-cyclyc and cyclic impulsive self-mappings

Consider a metric space $(X, d)$ and a composite self-mapping $T: \bigcup_{i\in \bar{p}} A_i \to \bigcup_{i\in \bar{p}} A_i$ of the form $T = T^+ \circ T^-$ (abbreviated in the following simply as $T = T^+ T^-$), where $A_i$ ; $i \in \bar{p}$ are $p(\geq 2)$ nonempty closed subsets of $X$ with $A_{np+i} \equiv A_i$ ; $\forall i \in \bar{p}$, $\forall n \in \mathbf{Z}_{0+}$ of distance $D = d(A_i, A_j) \geq 0$ between any two distinct subsets $A_i$ and $A_j$ of $X$ ; $\forall i, j(\neq i) \in \bar{p}$. Some useful types of such composite self-mappings for applications together with some of their properties in metric spaces are studied in this paper according to the following set of definitions:

**Definition 2.1**. The composite self-mapping $T: \bigcup_{i\in \bar{p}} A_i \to \bigcup_{i\in \bar{p}} A_i$ is said to be $p$-semi-cyclic if the following conditions hold:

1) $T^-: \bigcup_{i\in \bar{p}} A_i \to \bigcup_{i\in \bar{p}} A_i$ is such that $T^- A_i \subseteq A_{i+1}$ ; $\forall i \in \bar{p}$

2) $T^+: \bigcup_{i\in \bar{p}} A_i \to \bigcup_{i\in \bar{p}} A_i$ is such $T^+ A_i \subseteq A_i \cup A_{i+1}$ ; $\forall i \in \bar{p}$

3) $0 \leq d(T^- x, T^- y) \leq K d(x, y) + (1 - K) D$ ; $\forall x \in A_i$, $\forall y \in A_{i+1}$, $\forall i \in \bar{p}$ for some real constant $K \in \mathbf{R}_{0+}$

4) $d(Tx, Ty) = d(T^+(T^- x), T^+(T^- y)) \leq m(T^- x, T^- y) d(T^- x, T^- y)$ ; $\forall x \in A_i$, $\forall y \in A_{i+1}$, $\forall i \in \bar{p}$ for some function $m: (\bigcup_{i\in \bar{p}} A_i) \times (\bigcup_{i\in \bar{p}} A_i) \to \mathbf{R}_{0+}$.  □

Composite, in general non-continuous, self-mappings satisfying Definition 2.1 will be referred to in the following simply as *semi-cyclic impulsive self- mappings* because of its usefulness in describing mixed nonexpansive/contractive properties together with impulsive effects. See the conditions 3 and 4 of Definition 2.1. Examples are given later on in the paper concerning their application to the stability impulsive differential and difference equations and dynamic systems. Thus, the composite self-mapping of Definition 2.1 is being referred to be a $p-$ semi-cyclic impulsive self-mapping. Such a self-mapping is more general than the so-called $p-$ cyclic self-mapping $T = T^- : \bigcup_{i \in \overline{p}} A_i \to \bigcup_{i \in \overline{p}} A_i$ ( so that $T^+$ is identity) with $TA_i \subseteq A_{i+1}$; $\forall i \in \overline{p}$ (thus, $TA_i \subseteq A_i$ is not allowed in cyclic self-mappings while it is allowed in semi-cyclic ones). Further particular specifications of $p-$ semi-cyclic impulsive self-mappings are given in the subsequent definitions:

**Definition 2.2**. The $p-$ semi-cyclic impulsive self-mapping $T : \bigcup_{i \in \overline{p}} A_i \to \bigcup_{i \in \overline{p}} A_i$ is said to be nonexpansive (respectively, contractive) if $K \in [0,1]$ (respectively, if $K \in [0,1)$) and $m(T^- x, T^- y) \leq 1$; $\forall x \in A_i$, $\forall y \in A_{i+1}$, $\forall i \in \overline{p}$. □

Note that if $m(T^- x, T^- y) \leq 1$; $\forall x \in A_i$, $\forall y \in A_{i+1}$, $\forall i \in \overline{p}$ then

$$0 \leq d(Tx, Ty) \leq m(T^- x, T^- y) d(T^- x, T^- y) \leq m(T^- x, T^- y)(K d(x,y) + (1-K)D) \leq K d(x,y) + (1-K)D$$

$\forall x \in A_i$, $\forall y \in A_{i+1}$, $\forall i \in \overline{p}$. Note that a necessary condition for $T : \bigcup_{i \in \overline{p}} A_i \to \bigcup_{i \in \overline{p}} A_i$ to be a nonexpansive (respectively, contractive) $p-$ semi-cyclic impulsive self-mapping is that $T^- : \bigcup_{i \in \overline{p}} A_i \to \bigcup_{i \in \overline{p}} A_i$ be a nonexpansive (respectively, contractive) $p-$ cyclic impulsive self-mapping. Definitions 2.1 and 2.2 do not guarantee that $T : \bigcup_{i \in \overline{p}} A_i \to \bigcup_{i \in \overline{p}} A_i$ is $p-$ cyclic although $T^- : \bigcup_{i \in \overline{p}} A_i \to \bigcup_{i \in \overline{p}} A_i$ is always cyclic by virtue of the condition 1 of Definition 2.1. Note that the condition 4 of Definition 2.1 does not require $d(Tx, Ty) \geq D$; $\forall x \in A_i$, $\forall y \in A_{i+1}$, $\forall i \in \overline{p}$ and, from the conditions 1-2 of Definition 2.1 $TA_i \subseteq A_i \cup A_{i+1}$.

**Definition 2.3**. The $p-$ semi-cyclic impulsive self-mapping $T : \bigcup_{i \in \overline{p}} A_i \to \bigcup_{i \in \overline{p}} A_i$ is said to be $p-$ cyclic if it is $p-$ semi-cyclic and, furthermore, $d(Tx, Ty) \geq D$; $\forall x \in A_i$, $\forall y \in A_{i+1}$, $\forall i \in \overline{p}$. □

**Definition 2.4**. The $p-$ cyclic impulsive self-mapping $T : \bigcup_{i \in \overline{p}} A_i \to \bigcup_{i \in \overline{p}} A_i$ is said to be nonexpansive (respectively, contractive) $p-$ cyclic if $K \in [0,1]$ (respectively, if $K \in [0,1)$), $m(T^- x, T^- y) \leq 1$ and $d(Tx, Ty) \geq D$; $\forall x \in A_i$, $\forall y \in A_{i+1}$, $\forall i \in \overline{p}$. □

Note that if $D=0$, a nonexpansive $p$-semi-cyclic impulsive self-mapping may be also be simultaneously $p$-cyclic if for $x \in A_i$, $Tx \in A_{i+1}$ (instead of $Tx \in A_i \cup A_{i+1}$); $\forall i \in \overline{p}$. If the condition 4 is modified as follows

$$d(Tx, Ty) \leq Km(T^-x, T^-y)d(x,y) + (1 - Km(T^-x, T^-y))D \tag{2.1}$$

; $\forall x \in A_i$, $\forall y \in A_{i+1}$, $\forall i \in \overline{p}$, then $T: \bigcup_{i \in \overline{p}} A_i \to \bigcup_{i \in \overline{p}} A_i$ is said to be:

**Definition 2.5**. $T: \bigcup_{i \in \overline{p}} A_i \to \bigcup_{i \in \overline{p}} A_i$ is strictly $p$-semi-cyclic if (2.1) holds under the remaining conditions as in Definition 2.1. □

**Definition 2.6**. $T: \bigcup_{i \in \overline{p}} A_i \to \bigcup_{i \in \overline{p}} A_i$ is strictly nonexpansive (respectively, contractive) $p$-semi-cyclic if (2.1) holds under the remaining conditions as in Definition 2.2. □

**Definition 2.7**. $T: \bigcup_{i \in \overline{p}} A_i \to \bigcup_{i \in \overline{p}} A_i$ is strictly $p$-cyclic if (2.1) holds under the remaining conditions as in Definition 2.3. □

**Definition 2.8**. $T: \bigcup_{i \in \overline{p}} A_i \to \bigcup_{i \in \overline{p}} A_i$ is strictly nonexpansive (respectively, contractive) $p$-cyclic if (2.1) holds under the remaining conditions as in Definition 2.4. □

*Remark 2.9*. Note that if $m(T^-x, T^-y) \leq 1$; $\forall x \in A_i$, $\forall y \in A_{i+1}$, $\forall i \in \overline{p}$ then $m(T^-x, T^-y)(1-K)D \leq (1 - Km(T^-x, T^-y))D$; $\forall x \in A_i$, $\forall y \in A_{i+1}$, $\forall i \in \overline{p}$ and this holds if $D=0$ irrespective $m(T^-x, T^-y)$; $\forall x \in A_i$, $\forall y \in A_{i+1}$, $\forall i \in \overline{p}$. □

Thus, one has the following result from Remark 2.9:

**Proposition 2.10**. Assume that either $\bigcap_{i \in \overline{p}} A_i \neq \emptyset$ or $\bigcap_{i \in \overline{p}} A_i = \emptyset$ and $0 \leq m(T^-x, T^-y) \leq 1$; $\forall x \in A_i$, $\forall y \in A_{i+1}$, $\forall i \in \overline{p}$. Then, the self-mapping $T: \bigcup_{i \in \overline{p}} A_i \to \bigcup_{i \in \overline{p}} A_i$ is:

**(i)** Strictly $p$-semi-cyclic if it is $p$-semi-cyclic.

**(ii)** Strictly nonexpansive (respectively, contractive) $p$-semi-cyclic if it is nonexpansive (respectively, contractive) $p$-semi-cyclic.

**(iii)** Strictly $p$-cyclic if it is $p$-cyclic.

**(iv)** Strictly nonexpansive (respectively, contractive) $p$-cyclic if it is nonexpansive (respectively, contractive) $p$-cyclic. □

It is of interest the study of weaker properties than the above ones in an asymptotic context to be then able to investigate the asymptotic properties of distances for sequences $\{x_n\}_{n \in \mathbb{Z}_{0+}}$ of iterates built through $T: \bigcup_{i \in \bar{p}} A_i \to \bigcup_{i \in \bar{p}} A_i$ according to $x_{n+1} = T x_n$ for all $n \in \mathbf{Z}_{0+}$ and some $x_0 \in \bigcup_{i \in \bar{p}} A_i$ as well as the existence and uniqueness of fixed and best proximity points.

**Lemma 2.11.** Consider the $p$-semi-cyclic impulsive self-mapping $T: \bigcup_{i \in \bar{p}} A_i \to \bigcup_{i \in \bar{p}} A_i$ with $K \in [0,1]$ and define

$$m'(T^- x, T^- y) = m(T^- x, T^- y) - 1 \ , \ \delta_k(x) = m'(T^{(k+1)-} x, T^{k-} x)(K d(T^k x, T^{k-1} x) + (1-K)D)$$

for $x$ and $y$ in adjacent subsets $A_i$ and $A_{i+1}$ of $X$ for any $i \in \bar{p}$. Then, the following properties hold:

**(i)** The sequence $\{d(T^{k+np+j} x, T^{k+np+j-1} x)\}_{k \in \mathbf{Z}_{0+}}$ is bounded; $\forall k \in \mathbf{Z}_{0+}$, $\forall n \in \mathbf{Z}_+$, $\forall j \in \overline{p-1} \cup \{0\}$ if

$$-d(T^{k+1} x, T^k x) \le \sum_{i \in S_+(k,n,j)} \delta_{k+j+np-i}(x) - \sum_{i \in S_-(k,n,j)} \delta_{k+j+np-i}(x) < \infty \tag{2.2}$$

$\forall k \in \mathbf{Z}_{0+}$, $\forall n \in \mathbf{Z}_+$, $\forall j \in \overline{p-1} \cup \{0\}$, where

$$\delta_k(x) = m'(T^{(k+1)-} x, T^{k-} x)(K d(T^k x, T^{k-1} x) + (1-K)D) \ ; \ \forall k \in \mathbf{Z}_{0+} \tag{2.3}$$

$$S_+(k,n,j) = \{i \in \mathbf{Z}_+ : (i \le np + j) \wedge (m'(T^{(k+np+j-i+1)-} x, T^{(k+np+j-i)-} x)) > 0\}$$
$$S_-(k,n,j) = \{i \in \mathbf{Z}_+ : (i \le np + j) \wedge (-1 \le m'(T^{(k+np+j-i+1)-} x, T^{(k+np+j-i)-} x)) < 0\} \tag{2.4}$$

; $\forall k \in \mathbf{Z}_{0+}$, $\forall n \in \mathbf{Z}_+$, $\forall j \in \overline{p-1} \cup \{0\}$.

If, furthermore, $T: \bigcup_{i \in \bar{p}} A_i \to \bigcup_{i \in \bar{p}} A_i$ is if $p$-cyclic then the lower-bound in (2.2) is replaced with $D - d(T^{k+1} x, T^k x)$. If $T: \bigcup_{i \in \bar{p}} A_i \to \bigcup_{i \in \bar{p}} A_i$ is a nonexpansive $p$-semi-cyclic impulsive self-mapping (in particular, $p$-cyclic) then $\{d(T^{k+np+j} x, T^{k+np+j-1} x)\}_{k \in \mathbf{Z}_{0+}}$ is bounded; $\forall k \in \mathbf{Z}_{0+}$, $\forall n \in \mathbf{Z}_+$, $\forall j \in \overline{p-1} \cup \{0\}$.

**(ii)** If, furthermore, If $K \in [0,1)$ then

$$0 \le \limsup_{n \to \infty} d(T^{k+np+j} x, T^{k+np+j-1} x) \le D + \limsup_{n \to \infty} \sum_{i=1}^{np+j} \delta_{k+j+np-i}(x)$$
$$\le D + \limsup_{n \to \infty} \left( \sum_{i \in S_+(k,n,j)} \delta_{k+j+np-i}(x) - \sum_{i \in S_-(k,n,j)} \delta_{k+j+np-i}(x) \right) < \infty \tag{2.5}$$

; $\forall x \in \bigcup_{i \in \bar{p}} A_i$, $\forall j \in \overline{p-1} \cup \{0\}$. If, furthermore, $T: \bigcup_{i \in \bar{p}} A_i \to \bigcup_{i \in \bar{p}} A_i$ is if $p$-cyclic then the lower-bound in (2.5) is replaced with $D$.

If $T: \bigcup_{i \in \bar{p}} A_i \to \bigcup_{i \in \bar{p}} A_i$ is contractive $p$-semi-cyclic then $0 \le \limsup_{n \to \infty} d(T^{k+np+j} x, T^{k+np+j-1} x) \le D$;

$\forall k \in \mathbf{Z}_{0+}$, $\forall j \in \overline{p-1} \cup \{0\}$, $\forall x \in \bigcup_{i \in \bar{p}} A_i$.

If $T: \bigcup_{i \in \bar{p}} A_i \to \bigcup_{i \in \bar{p}} A_i$ is contractive $p$-cyclic then there exists $\lim_{n \to \infty} d(T^{k+np+1} x, T^{k+np} x) = D$;

$\forall x \in \bigcup_{i \in \bar{p}} A_i$.

*Proof*: Build a sequence of iterates $\{T^k x\}_{k \in \mathbf{Z}}$ according to $TT^{k-1}x = T^+T^-T^{k-1}x$ with $T^{0-}x = x$, $T^0 x = T^{0+}T^{0-}x = x$, for any given $x \in A_i$ and any $i \in \bar{p}$ that is $T = T^{0+} = T^{0-} = id$ so that

$$d(T^{k+1}x, T^k x) \leq (1 + m'(T^{(k+1)-}x, T^{k-}x)) d(T^{(k+1)-}x, T^{k-}x)$$
$$\leq (1 + m'(T^{(k+1)-}x, T^{k-}x))(Kd(T^k x, T^{k-1}x) + (1-K)D)$$
$$= Kd(T^k x, T^{k-1}x) + (1-K)D + \delta_k(x) \qquad (2.6)$$

; $\forall k \in \mathbf{Z}_{0+}$. Through a recursive calculation with (2.2), one gets:

$$0 \leq d(T^{k+np+j}x, T^{k+np+j-1}x) \leq Kd(T^{k+np+j-1}x, T^{k+np+j-2}x) + (1-K)D + \delta_{k+np+j-1}(x)$$
$$\leq K^2 d(T^{k+np+j-2}x, T^{k+np+j-3}x) + K[(1-K)D + \delta_{k+np+j-2}(x)] + (1-K)D + \delta_{k+np+j-1}(x)$$
$$\leq \ldots \leq K^{np+j-1} d(T^{k+1}x, T^k x) + (1-K^{np+j-1})D + \sum_{i=1}^{np+j} K^i \delta_{k+np+j-i}(x) \qquad (2.7)$$

; $\forall k \in \mathbf{Z}_{0+}$, $\forall n \in \mathbf{Z}_+$, $\forall j \in \overline{p-1} \cup \{0\}$. If $K = 1$ then

$$0 \leq d(T^{k+np+j}x, T^{k+np+j-1}x) \leq d(T^{k+1}x, T^k x) + \sum_{i \in S_+(k,n,j)} \delta_{k+np+j-i}(x) - \sum_{i \in S_-(k,n,j)} \delta_{k+np+j-i}(x)$$
$$(2.8)$$

; $\forall k \in \mathbf{Z}_{0+}$, $\forall n \in \mathbf{Z}_+$, $\forall j \in \overline{p-1} \cup \{0\}$. Take any $k \in \mathbf{Z}_{0+}$, any $n \in \mathbf{Z}_+$ and any $x \in \bigcup_{i \in \bar{p}} A_i$. Since $d(T^{k+1}x, T^k x)$ is finite and (2.2) holds it follows that $0 \leq d(T^{k+np+j+1}x, T^{k+np+j}x) < \infty$. If, in addition, $T : \bigcup_{i \in \bar{p}} A_i \to \bigcup_{i \in \bar{p}} A_i$ is $p-$ cyclic then the zero lower-bound of (2.5) is replaced with $D$. If $T : \bigcup_{i \in \bar{p}} A_i \to \bigcup_{i \in \bar{p}} A_i$ is $p-$ semi-cyclic (in particular, $p-$ cyclic) nonexpansive then (2.2) always holds since $m(T^{(k+np+j+i)-}x, T^{(k+np+j+i-1)-}x) \leq 1$, $-1 \leq m'(T^{(k+j+np-i+1)-}x, T^{(k+j+np-i)-}x) \leq 0$ so that

$$\sum_{i \in S_+(k,n,j)} \delta_{k+j+np-i}(x) - \sum_{i \in S_-(k,n,j)} \delta_{k+j+np-i}(x) = \sum_{i \in S_-(k,n,j)} \delta_{k+j+np-i}(x) \leq 0 \qquad (2.9)$$

if $m(T^{(k+j+np-i+1)-}x, T^{(k+j+np-i)-}x) = 1$ and $\{d(T^{k+np+j+1}x, T^{k+np+j}x)\}_{k \in \mathbf{Z}_{0+}}$ is always bounded; $\forall k \in \mathbf{Z}_{0+}$, $\forall n \in \mathbf{Z}_+$, $\forall j \in \overline{p-1} \cup \{0\}$. Property (i) has been proven. If $K \in [0,1)$ then

$$0 \leq d(T^{k+np+j}x, T^{k+np+j-1}x) \leq K^{np+j-1} d(T^{k+1}x, T^k x) + (1-K^{np+j-1})D$$
$$+ \sum_{i \in S_+(k,n,j)} \delta_{k+j+np-i}(x) - \sum_{i \in S_-(k,n,j)} \delta_{k+j+np-i}(x) \qquad (2.10)$$

$$0 \leq \limsup_{n \to \infty} d(T^{k+np+j}x, T^{k+np+j-1}x) \leq D + \limsup_{n \to \infty} \sum_{i=1}^{np+j} \delta_{k+j+np-i}(x) \qquad (2.11)$$

If, in addition, $T : \bigcup_{i \in \bar{p}} A_i \to \bigcup_{i \in \bar{p}} A_i$ is $p-$ cyclic then the zero lower-bound of (2.9)-(2.11) is replaced with $D$.

If $T:\bigcup_{i\in\bar{p}} A_i \to \bigcup_{i\in\bar{p}} A_i$ is contractive $p$- semi-cyclic then (2.11) becomes

$0 \le \limsup\limits_{n\to\infty} d\left(T^{k+np+1}x, T^{k+np}x\right) \le D$ from (2.9). If, in addition, $T:\bigcup_{i\in\bar{p}} A_i \to \bigcup_{i\in\bar{p}} A_i$ is contractive $p-$ cyclic then $D \le \limsup\limits_{n\to\infty} d\left(T^{k+np+1}x, T^{k+np}x\right) \le D$; $\forall x \in \bigcup_{i\in\bar{p}} A_i$ so that there is

$\lim\limits_{n\to\infty} d\left(T^{k+np+1}x, T^{k+np}x\right) = D$; $\forall x \in \bigcup_{i\in\bar{p}} A_i$. Property (ii) has been proven. □

The following result establishes an asymptotic property of limits superiors of distances of consecutive iterated points which implies that $T:\bigcup_{i\in\bar{p}} A_i \to \bigcup_{i\in\bar{p}} A_i$ is asymptotically contractive the limit

$\lim\limits_{n\to\infty}\left(\sum_{k=0}^{np+j-2}\left(\prod_{\ell=k}^{np+j-2}[K_{\ell+i}]\right)\left(m\left(T^{(k+1)-}x, T^{k-}x\right)-1\right)\right) = 0$; $\forall x \in \bigcup_{i\in\bar{p}} A_i$, $\forall j \in \overline{p-1}\cup\{0\}$ exists.

In particular, it is not required that $m(x,y) \le 1$ for any $x \in A_i$, $y \in A_{i+1}$; $\forall i \in \bar{p}$ as in contractive and, in general, nonexpansive $p$- semi-cyclic impulsive self-mappings.

**Theorem 2.12**. Consider the following generalization of condition 3 of Definition 2.1:

$$D \le d\left(T^{2-}x, T^{-}y\right) \le K_i d(Tx, x) + (1-K_i)D \tag{2.12}$$

for any given $x \in A_i$; $\forall i \in \bar{p}$ and define $\bar{K} = \prod_{i=1}^{p-1}[K_i]$. Define

$$\hat{K} = \bar{K} \sup_{x\in\bigcup_{i\in\bar{p}} A_i} \max_{n\in\mathbb{Z}_{0+}} \left(\prod_{i=np+1}^{(n+1)p-1}\left[m\left(T^{(i+1)-}x, T^{i-}x\right)\right]\right)$$

$$= \prod_{i=1}^{p-1}[K_i] \sup_{x\in\bigcup_{k\in\bar{p}} A_k} \max_{n\in\mathbb{Z}_{0+}} \left(\prod_{i=np+1}^{(n+1)p-1}\left[m\left(T^{(i+1)-}x, T^{i-}x\right)\right]\right) \tag{2.13}$$

such that $\hat{K} \in [0,1)$. Then, the following properties hold:

**(i)** $D_0 \le \limsup\limits_{n\to\infty} d\left(T^{np+j}x, T^{np+j-1}x\right) \le \left(1+\dfrac{1}{1-\hat{K}}\left(\prod_{\ell=0}^{j-1}[K_{i+\ell}]\right)\sup_{x\in\bigcup_{k\in\bar{p}} A_k}\max_{\ell\in\mathbb{Z}_{0+}}\left|m'\left(T^{(\ell+1)-}x, T^{\ell-}x\right)\right|\right)D$

; $\forall x \in \bigcup_{i\in\bar{p}} A_i$, $\forall j \in \overline{p-1}\cup\{0\}$ (2.14)

$D_0 \le d\left(T^{np+j}x, T^{np+j-1}x\right) \le \left(\prod_{\ell=0}^{j-1}[K_{i+\ell}]\right)\hat{K}^n d(Tx, x)$

$+ \left[\left(1-\left(\prod_{\ell=0}^{j-1}[K_{i+\ell}]\right)\hat{K}^n\right) + \dfrac{1-\hat{K}^n}{1-\hat{K}}\left(\prod_{\ell=0}^{j-1}[K_{i+\ell}]\right)\sup_{x\in\bigcup_{k\in\bar{p}} A_k}\max_{\ell\in\mathbb{Z}_{0+}}\left|m'\left(T^{(\ell+1)-}x, T^{\ell-}x\right)\right|\right]D$

; $\forall x \in A_i$, $\forall i \in \bar{p}$, $\forall j \in \overline{p-1}\cup\{0\}$ (2.15)

where $D_0 = 0$ if $T:\bigcup_{i\in\bar{p}} A_i \to \bigcup_{i\in\bar{p}} A_i$ is $p-$ pre-cyclic and $D_0 = D$ if $T:\bigcup_{i\in\bar{p}} A_i \to \bigcup_{i\in\bar{p}} A_i$ is $p-$ cyclic.

**(ii)** If, furthermore, there is $\varepsilon_0 (\ge -1) \in \mathbf{R}$ such that

$$\limsup_{n\to\infty}\left(\sum_{k=0}^{np+j-2}\left(\prod_{\ell=k-1}^{np+j-2}[K_{\ell+i}]\right)\left(m\left(T^{(k+1)-}x,T^{k-}x\right)-1\right)\right)\le \varepsilon_0\,;\,\forall x\in\bigcup_{i\in\overline{p}}A_i\,,\,\forall j\in\overline{p-1}\cup\{0\} \quad (2.16)$$

then

$$D_0\le \limsup_{n\to\infty} d\left(T^{np+j}x,T^{np+j-1}x\right)\le D(1+\varepsilon_0)\,;\,\forall x\in\bigcup_{i\in\overline{p}}A_i,\,\forall j\in\overline{p-1}\cup\{0\} \quad (2.17)$$

*Proof*: Since $\hat{K}\in[0,1)$, one has through iterative calculation via (2.12) that

$$d(T^2x,Tx)\le m(T^{2-}x,T^-x)(K_i d(Tx,x)+(1-K)D)$$

$$=(m(T^{2-}x,T^-x)K_i)d(Tx,x)+(1-m(T^{2-}x,T^-x)K_i)D+m'(T^{2-}x,T^-x)D$$

$$\ldots\ldots\ldots\ldots\ldots\ldots\ldots\ldots\ldots\ldots\ldots\ldots\ldots\ldots\ldots\ldots\ldots$$

$$d(T^px,T^{p-1}x)\le\left(\prod_{i=1}^{p-1}\left[m\left(T^{(i+1)-}x,T^{i-}x\right)\right]\right)\overline{K}d(Tx,x)$$

$$+\left(1-\left(\prod_{i=1}^{p-1}\left[m\left(T^{(i+1)-}x,T^{i-}x\right)\right]\right)\overline{K}\right)D+\left(\sum_{k=0}^{p-2}\left(\prod_{\ell=k-1}^{p-2}[K_{\ell+i}]\right)m'\left(T^{(k+1)-}x,T^{k-}x\right)\right)D$$

$$\ldots\ldots\ldots\ldots\ldots\ldots\ldots\ldots\ldots\ldots\ldots\ldots\ldots\ldots\ldots\ldots\ldots$$

$$d\left(T^{np+j}x,T^{np+j-1}x\right)\le\left(\prod_{i=1}^{np+j-1}\left[m\left(T^{(i+1)-}x,T^{i-}x\right)\right]\right)\left(\prod_{\ell=0}^{j-1}[K_{i+\ell}]\right)\overline{K}^n d(Tx,x)$$

$$+\left(1-\left(\prod_{i=1}^{np+j-1}\left[m\left(T^{(i+1)-}x,T^{i-}x\right)\right]\right)\left(\prod_{\ell=0}^{j-1}[K_{i+\ell}]\right)\overline{K}^n\right)D+\left(\sum_{k=0}^{np+j-2}\left(\prod_{\ell=k-1}^{np+j-2}[K_{\ell+i}]\right)m'\left(T^{(k+1)-}x,T^{k-}x\right)\right)D$$

$$;\,\forall x\in A_i,\,\forall i\in\overline{p},\,\forall j\in\overline{p-1}\cup\{0\} \quad (2.18)$$

with the convention $\left(\prod_{\ell=0}^{-1}[K_{i+\ell}]\right)=1$; $\forall i\in\overline{p}$. Then, one gets (2.14)-(2.15) and Property (i) has been proven. To prove property (ii), use the indicator sets (2.4) and, since $m'(T^{2-}x,T^-x)\ge-1$; $\forall x\in\bigcup_{i\in\overline{p}}A_i$, one also gets from (2.12)-(2.13):

$$D_0+\left[\liminf_{n\to\infty}\left(\sum_{k\in S_-(k,n,j-2)}\left(\prod_{\ell=k}^{np+j-2}[K_{\ell+i}]\right)\left|m'\left(T^{(k+1)-}x,T^{k-}x\right)\right|\right)\right]D$$

$$\le\limsup_{n\to\infty} d\left(T^{np+j}x,T^{np+j-1}x\right)+\left[\liminf_{n\to\infty}\left(\sum_{k\in S_-(k,n,j-2)}\left(\prod_{\ell=k-1}^{np+j-2}[K_{\ell+i}]\right)\left|m'\left(T^{(k+1)-}x,T^{k-}x\right)\right|\right)\right]D$$

$$\le D\left[\limsup_{n\to\infty}\left(1+\left(\sum_{k\in S_+(k,n,j-2)}\left(\prod_{\ell=k-1}^{np+j-2}[K_{\ell+i}]\right)m'\left(T^{(k+1)-}x,T^{k-}x\right)\right)\right)\right]$$

$$;\,\forall x\in\bigcup_{i\in\overline{p}}A_i,\,\forall j\in\overline{p-1}\cup\{0\} \quad (2.19)$$

and (2.17), and then Property (ii), follows from (2.16). □

Note from Theorem 2.12, Eq. (2.17) that if $D_0=D=0$, that is $\bigcap_{i\in\overline{p}}A_i\ne\varnothing$, and $\varepsilon_0\in[-1,\infty)$ then $\exists\lim_{n\to\infty} d\left(T^{np+j}x,T^{np+j-1}x\right)=0$; $\forall x\in\bigcup_{i\in\overline{p}}A_i,\,\forall j\in\overline{p-1}\cup\{0\}$ from (2.17) since $\hat{K}\in[0,1)$. In this case, $T: \bigcup_{i \in \bar{p}} A_i \to \bigcup_{i \in \bar{p}} A_i$ is an asymptotically contractive $p$-cyclic (and also $p$-semi-cyclic since $D = 0$) self-mapping on the union on intersecting closed subsets of $X$. A close property follows if $D_0 = D \neq 0$ and $\varepsilon_0 = 0$ implying from (2.17) that

$$\limsup_{n \to \infty} \left( \sum_{k=0}^{np+j-2} \left( \prod_{\ell=k-1}^{np+j-2} [K_{\ell+i}] \right) \left( m\left(T^{(k+1)-}x, T^{k-}x\right) - 1 \right) \right)$$
$$= \lim_{n \to \infty} \left( \sum_{k=0}^{np+j-2} \left( \prod_{\ell=k-1}^{np+j-2} [K_{\ell+i}] \right) \left( m\left(T^{(k+1)-}x, T^{k-}x\right) - 1 \right) \right) = 0 \,;\, \forall x \in \bigcup_{i \in \bar{p}} A_i \,,\, \forall j \in \overline{p-1} \cup \{0\}$$

and leading to $\exists \lim_{n \to \infty} d\left(T^{np+j}x, T^{np+j-1}x\right) = D$ such that $T: \bigcup_{i \in \bar{p}} A_i \to \bigcup_{i \in \bar{p}} A_i$ is a contractive $p$-cyclic self-mapping on the union on disjoint closed subsets of $X$. The above discussion is summarized in the subsequent result:

**Corollary 2.13.** Assume that (2.12) holds with $\hat{K}$ defined in (2.13) being in $[0, 1)$ and assume also that:

$$\infty > \varepsilon_0 \geq \max \left( \limsup_{n \to \infty} \left( \sum_{k=0}^{np+j-2} \left( \prod_{\ell=k-1}^{np+j-2} [K_{\ell+i}] \right) \left( m\left(T^{(k+1)-}x, T^{k-}x\right) - 1 \right), -1 \right) \right)$$

$;\, \forall x \in \bigcup_{i \in \bar{p}} A_i \,,\, \forall j \in \overline{p-1} \cup \{0\}$. Then, the following properties hold:

**(i)** If $\bigcap_{i \in \bar{p}} A_i \neq \emptyset$ then $T: \bigcup_{i \in \bar{p}} A_i \to \bigcup_{i \in \bar{p}} A_i$ is an asymptotically contractive $p$-cyclic impulsive self-mapping so that there is the limit
$$\lim_{n \to \infty} d\left(T^{np+j}x, T^{np+j-1}x\right) = 0 \,;\, \forall x \in \bigcup_{i \in \bar{p}} A_i \,,\, \forall j \in \overline{p-1} \cup \{0\} \,.$$

**(ii)** If $\bigcap_{i \in \bar{p}} A_i = \emptyset$, $d(Tx, Ty) \geq D$; $\forall x \in A_i$, $\forall y \in A_{i+1}$, $\forall i \in \bar{p}$ and the following limit exists:

$$\lim_{n \to \infty} \left( \sum_{k=0}^{np+j-2} \left( \prod_{\ell=k-1}^{np+j-2} [K_{\ell+i}] \right) \left( m\left(T^{(k+1)-}x, T^{k-}x\right) - 1 \right) \right) = 0 \,;\, \forall x \in \bigcup_{i \in \bar{p}} A_i \,,\, \forall j \in \overline{p-1} \cup \{0\} \quad (2.20)$$

then $T: \bigcup_{i \in \bar{p}} A_i \to \bigcup_{i \in \bar{p}} A_i$ is an asymptotically contractive $p$-cyclic impulsive self-mapping so that the limit
$$\lim_{n \to \infty} d\left(T^{np+j}x, T^{np+j-1}x\right) = D \,;\, \forall x \in \bigcup_{i \in \bar{p}} A_i \,,\, \forall j \in \overline{p-1} \cup \{0\} \text{ exists.} \qquad \square$$

A particular result obtained from Theorem 2.12 follows for contractive $p$-semi-cyclic and $p$-cyclic impulsive self-mappings $T: \bigcup_{i \in \bar{p}} A_i \to \bigcup_{i \in \bar{p}} A_i$:

**Corollary 2.14.** Theorem 2.12 holds with $D_0 = 0$ if $T: \bigcup_{i \in \bar{p}} A_i \to \bigcup_{i \in \bar{p}} A_i$ is contractive $p$-semi-cyclic and with $D_0 = D$ if the impulsive self-mapping $T: \bigcup_{i \in \bar{p}} A_i \to \bigcup_{i \in \bar{p}} A_i$ is contractive $p$-cyclic provided that $\overline{K} = \prod_{i=1}^{p-1} [K_i] \in [0, 1)$.

*Proof*: It is a direct consequence of Theorem 2.12 since $\overline{K} = \prod_{i=1}^{p-1}[K_i] \in [0,1)$ implies that $\hat{K} \in [0,1)$ since $m(T^-x, T^-y) \leq 1$; $\forall x \in A_i$, $\forall y \in A_{i+1}$, $\forall i \in \overline{p}$. □

**Remark 2.15**. Note that if $T: \bigcup_{i \in \overline{p}} A_i \to \bigcup_{i \in \overline{p}} A_i$ is a nonexpansive $p$-cyclic impulsive self-mapping, the following constraints hold:

$$m(Tx^-, T y^-) \leq 1 \; ; \; D \leq m(Tx^-, T y^-)(K d(x, y) - D) + m(Tx^-, T y^-)D \; ; \; \forall x \in A_i, \; \forall y \in A_{i+1}, \; \forall i \in \overline{p}$$

, and equivalently,

$$1 \geq m(Tx^-, T y^-) \geq \frac{D}{Kd(x, y) + (1 - K)D} = \frac{D}{D + K(d(x, y) - D)} \tag{2.21}$$

implying that

a) $1 \geq m(Tx^-, T y^-) \geq 0$; $\forall x \in A_i$, $\forall y \in A_{i+1}$, $\forall i \in \overline{p}$ if $D = 0$, i.e. if the sets $A_i$ intersect; $\forall i \in \overline{p}$

b) $m(Tx^-, T y^-) = 1$ if $d(x, y) = D$; i.e. for best proximity points associated with any two adjacent disjoint subsets $A_i$, $y \in A_{i+1}$ for $i \in \overline{p}$. □

Some functions are now defined to evaluate the nonexpansive and contractive properties of the impulsive self-mapping $T: \bigcup_{i \in \overline{p}} A_i \to \bigcup_{i \in \overline{p}} A_i$ which take into account the most general case that the constants in Definition 2.1 (3-4) can be generalized to be set dependent and point-dependent leading to a combined extended constraint as follows:

$$d(T^2 x, Tx) \leq K_i(x, Tx) m(T^{2-}x, T^-x) d(x, Tx) + m(T^{2-}x, T^-x)(1 - K_i(x, Tx))D \; ; \; \forall x \in \bigcup_{i \in \overline{p}} A_i \tag{2.22}$$

so that

$$\hat{K}^{(j)}(x, Tx) = \left( \prod_{i=1}^{p-1} \left[ m\left(T^{(i+jp+1)-}x, T^{(i+jp)-}x\right) K_i\left(T^{i+jp}x, T^{i+jp-1}x\right) \right] \right) \hat{K}^{(j-1)}(x) \; ; \; \forall x \in \bigcup_{i \in \overline{p}} A_i, \; \forall j \in \mathbf{Z}_+ \tag{2.23}$$

with $x = T^0 x$ and initial, in general, point-dependent value

$$\hat{K}^{(0)}(x, Tx) = \prod_{i=1}^{p-1} \left[ m\left(T^{(i+1)-}x, T^{(i)-}x\right) K_i\left(T^i x, T^{i-1}x\right) \right] \; ; \; \forall x \in \bigcup_{i \in \overline{p}} A_i, \; \forall j \in \mathbf{Z}_+ \tag{2.24}$$

**3. Convergence of the iterations to best proximity points and fixed points**

Important results about convergence of iterated sequences of 2-cyclic self-mappings to unique best proximity points were firstly stated and proven in [11] and then widely used in the literature. Some of them are quoted here to be then used in the context of this paper. Consider a metric space $(X, d)$ with nonempty subsets $A, B \subset X$ such that $D = d(A, B) \geq 0$. The following basic results have been proven in the existing background literature:

*Result 1* [11]. Let $(X, d)$ be a metric space and let $A$ and $B$ be subsets of $X$. Then, if $A$ is compact and $B$ is approximatively compact with respect to $A$ (i.e. $d(y, x_n) \to d(y, B)$ as $n \to \infty$ for each sequence $\{x_n\}_{n \in \mathbf{Z}_{0+}} \subset B$ for some $y \in A$) then $A^o = \{x \in A : d(x, y') = D \text{ for some } y' \in B\}$ and $B^o = \{x \in B : d(x', y) = D \text{ for some } x' \in A\}$ are nonempty. □

It is known that if $A$ and $B$ are both compact then $A$ (respectively, $B$) is approximatively compact which respect to $B$ (respectively, with respect to $A$).

*Result 2* [11]. Let $(X, \|\ \|)$ be a reflexive Banach space and let $A$ be a nonempty, closed, bounded and convex subset of $X$ and let $B$ be a nonempty, closed and convex subset of $X$. Then, the sets of best proximity points $A^o$ and $B^o$ are nonempty. □

*Result 3* [11]. Let $(X, d)$ be a metric space, let $A$ and $B$ be nonempty closed subsets of $X$ and let $T : A \cup B \to A \cup B$ be a 2-cyclic contraction. If either $A$ is boundedly compact (that is if any bounded sequence $\{x_n\}_{n \in \mathbf{Z}_{0+}} \subset A$ has a subsequence converging to a point of $A$) or $B$ is boundedly compact, then there is $x \in A \cup B$ such that $d(x, Tx) = D$. □

**Remark 3.1.** It is known that if $A \subset X$ is boundedly compact then it is approximatively compact. Also, a closed set $A$ of a normed space is boundedly compact if it is locally compact (the inverse is not true in separable Hilbert spaces), equivalently, if and only if the closure of each bounded subset $C \subset A$ is compact and contained in $A$. If $(X, d)$ is a linear metric space, a closed subset $A \subset X$ is boundedly compact if each bounded $C \subset A$ is relatively compact. It turns out that if $A \subset X$ is closed and bounded then it is relatively compact and approximatively compact with respect to itself. It also turns out that if $(X, d)$ is a complete metric space and the metric is homogeneous and translation-invariant then $(X, d)$ is a linear metric space and $(X, \|\ \|)$ is also a Banach space with $\|\ \|$ being the norm induced by the metric $d$ (note that, since the metric is homogeneous and translation–invariant, since $(X, d)$ is a linear metric space, it induces a norm). In such a Banach space, if $A \subset X$ is bounded and closed, then $A$ is boundedly compact and then approximatively compact with respect to itself.

*Result 4* [11]. Let $(X, \|\ \|)$ be a uniformly convex Banach space, let $A$ be a nonempty closed and convex subset of $X$ and let $B$ be a nonempty closed subset of $X$. Let sequences $\{x_n\}_{n \in \mathbf{Z}_{0+}} \subset A$, $\{z_n\}_{n \in \mathbf{Z}_{0+}} \subset A$ and $\{y_n\}_{n \in \mathbf{Z}_{0+}} \subset B$ satisfy $\|x_n - y_n\| \to D$ and $\|z_n - y_n\| \to D$ as $n \to \infty$. Then $\|z_n - x_n\| \to 0$ as $n \to \infty$. □

It is known that a uniformly convex Banach space $(X, \|\;\|)$ is reflexive and that a Banach space is a complete metric space $(X, d)$ with respect to the norm-induced distance.

*Result 5* [11]. If $(X, d)$ is a complete metric space, $T: A \cup B \to A \cup B$ be a $2$-cyclic contraction, where $A$ and $B$ are nonempty closed subsets of $X$, and the sequence $\{x_n\}_{n \in \mathbf{Z}_{0+}}$ generated as $x_{n+1} = Tx_n$; $\forall n \in \mathbf{Z}_+$ for a given $x_0 \in A$, has a convergent subsequence $\{x_{2n_k}\}_{n_k \in \mathbf{Z}_{0+}} \subset \{x_{2n}\}_{n \in \mathbf{Z}_{0+}} \subset \{x_n\}_{n \in \mathbf{Z}_{0+}}$ in $A$ then there is $x \in A \cup B$ such that $d(x, Tx) = D$. □

Sufficiency-type results follow below concerning the convergence of iterated sequences being generated by contractive and strictly contractive $p$- semi-cyclic self-mappings, which are asymptotically $p$-cyclic, to best proximity or fixed points:

**Theorem 3.2**. Assume that $(X, \|\;\|)$ is a uniformly convex Banach space so that $(X, d)$ is also a complete metric space if $d: X \times X \to \mathbf{R}_{0+}$ is the norm-induced metric. Assume, in addition, that $T: \bigcup_{i \in \bar{p}} A_i \to \bigcup_{i \in \bar{p}} A_i$ is a $p$- semi-cyclic impulsive self-mapping, where $A_i \subset X$; $\forall i \in \bar{p}$ are nonempty, closed and convex subsets of $X$ and assume also that

(1) Either the constraint (2.22) holds subject to (2.23)-(2.24) provided that the limit $\lim_{n \to \infty} \hat{K}^{(n)}(x, Tx) = 0$; $\forall x \in \bigcup_{i \in \bar{p}} A_i$ exists and $m: \left(\bigcup_{i \in \bar{p}} A_i\right) \times \left(\bigcup_{i \in \bar{p}} A_i\right) \to \mathbf{R}_{0+}$ satisfies the identity:

$$m\left(T^{((n+1)p-1)-}x, T^{((n+1)p-2)-}x\right) = 1 + \varepsilon_n - \frac{1}{\hat{K}^{(1)}(x, Tx)\left(\sum_{k=np}^{(n+1)p-2}\left(\prod_{\ell=(n+1)p+i-3}^{(n+1)p-2}\left[K_\ell\left(T^{((n+1)p-3)}x, T^{(n+1)p-2}x\right)\right]\right)\right)}$$

$$\times \left(\hat{K}^{(1)}(x, Tx)\left(\sum_{k=np}^{(n+1)p-3}\left(\prod_{\ell=k+i-1}^{(n+1)p-2}\left[K_\ell\left(T^{(k-1)}x, T^k x\right)\right]\right)\right)\left(m\left(T^{(k+1)-}x, T^{k-}x\right) - 1\right)\right)$$

$$+ \sum_{j=0}^{n-1} \hat{K}^{(n-j)}(x, Tx)\left(\sum_{k=jp}^{(j+1)p-2}\left(\prod_{\ell=k+i-1}^{(j+1)p-2}\left[K_\ell\left(T^{(k-1)}x, T^k x\right)\right]\right)\left(m\left(T^{(k+1)-}x, T^{k-}x\right) - 1\right)\right)$$

(2) For each given $x \in A_i$ for any $i \in \bar{p}$, there is a finite $k_i = k_i(x) \in \mathbf{Z}_{0+}$ such that $\liminf_{n \to \infty} T^{np+k_i(x)} \in A_{i+1}$ (i.e. the $p$- semi-cyclic impulsive self-mapping $T: \bigcup_{i \in \bar{p}} A_i \to \bigcup_{i \in \bar{p}} A_i$ is also an asymptotically $p$- cyclic one).

Then, $T: \bigcup_{i \in \bar{p}} A_i \to \bigcup_{i \in \bar{p}} A_i$ is either an asymptotically contractive or a strictly contractive $p$- semi-cyclic impulsive self-mapping and, furthermore, the following properties hold:

**(i)** The limits below exist:

$$\lim_{n \to \infty} d\left(T^{(n+1)p}x, T^{(n+j)p+j}x\right) = D; \; \forall x \in A_i, \forall j \in \bar{k}_i, \forall i \in \bar{p} \tag{3.1}$$

$$\lim_{n\to\infty} d\left(T^{(n+1)p+\bar{k}_i+1}x, T^{(n+j)p+\bar{k}_i}x\right) = 0 \; ; \; \forall x \in A_i, \; \forall i \in \bar{p} \qquad (3.2)$$

where $\bar{k}_i = \sup_{x \in A_i} k_i(x); \; \forall i \in \bar{p}$. Furthermore, $\{T^{np}x\}_{n \in \mathbf{Z}_+} \to z_i$, $\{T^{np+j}x\}_{n \in \mathbf{Z}_+} \to Tz_i^{(j)}$ for any given $x \in A_i$ with $\{T^{np+j}x\}_{n \in \mathbf{Z}_+} \subset A_i \cup A_{i+1}; \; \forall j \in \overline{\bar{k}_i}$, $\lim_{n\to\infty} T^{np+\bar{k}_i}x \subset A_{i+1}$, $z_i \in A_i$, $z_i^{(j)} \in A_i$; $\forall j \in \overline{\bar{k}_i} - 1$, and $z_{i+1} = Tz_i^{(k_i)} \in A_{i+1}$; $\forall i \in \bar{p}$. The points $z_i$ and $z_{i+1}$ are unique best proximity points in $A_i$ and $A_{i+1}$; $\forall \overline{ip}$ of $T : \bigcup_{i \in \bar{p}} A_i \to \bigcup_{i \in \bar{p}} A_i$ and there is a unique limiting set

$$\left( z_1, z_1^{(1)} = Tz_1, \ldots, z_2 = z_1^{(\bar{k}_1)} = T^{\bar{k}_1}z_1, \ldots, z_p, z_p^{(1)} = Tz_p, \ldots, z_p^{(\bar{k}_p - 1)} = T^{\bar{k}_p - 1}z_p \right) \subset A_1^{\bar{k}_1} \times \ldots \times A_1^{\bar{k}_p} \qquad (3.3)$$

If $\bigcap_{i \in \bar{p}} A_i \neq \varnothing$ then the $p$ best proximity points $z_i = z \in \bigcap_{j \in \bar{p}} A_j$; $\forall i \in \bar{p}$ become a unique fixed point $z$ of $T : \bigcup_{i \in \bar{p}} A_i \to \bigcup_{i \in \bar{p}} A_i$.

 (ii) Assume that the constraint (2.12) holds, subject to either (2.20), or to (2.22), with $\bar{K} = \prod_{i=1}^{p-1}[K_i]$ and $\hat{K} \in [0, 1)$ defined in (2.13). Assume that in addition, that for each $x \in A_i$ for any $i \in \bar{p}$, it exists a finite $k_i = k_i(x) \in \mathbf{Z}_{0+}$ such that $\liminf_{n \to \infty} T^{np+k_i(x)} \in A_{i+1}$ with $\bar{k}_i = \sup_{x \in A_i} k_i(x); \; \forall i \in \bar{p}$. Then, Property (i) still holds.

**Remarks 3.3**.
(1) Note that if the self-mapping $T : \bigcup_{i \in \bar{p}} A_i \to \bigcup_{i \in \bar{p}} A_i$ is an asymptotic $p-$ cyclic impulsive one then the limiting set (3.3) of Theorem 3.2 can only contain points which are not best proximity points in bounded subsets $A_i$ of $X$ whose diameter is not smaller than $D$.

(2) Under the conditions of Theorem 3.2, if $T : \bigcup_{i \in \bar{p}} A_i \to \bigcup_{i \in \bar{p}} A_i$ is, in particular, a contractive or strictly contractive $p-$ cyclic impulsive self-mapping, then the limiting set (3.3) only contains best proximity points, that is, it is of the form $(z_1,, z_2, \ldots, z_p)$. If $\bigcap_{i \in \bar{p}} A_i \neq \varnothing$ then such a set reduces to a unique best proximity point $z \in \bigcap_{i \in \bar{p}} A_i$.

(3) Note that Theorem 3.2 can be formulated also for a complete metric space $(X, d)$ with an homogeneous translation-invariant metric $d : X \times X \to \mathbf{R}_{0+}$ being equivalent to a Banach space $(X, \|\; \|)$, where $\|\; \|$ is the metric-induced norm, which is uniformly convex so that is also a complete. Note that such a statement is well-posed since a norm-induced metric exists if such a metric is homogeneous and translation-invariant. □

It turns out that Theorem 3.2 and Remarks 3.3 also hold if $T : \bigcup_{i \in \bar{p}} A_i \to \bigcup_{i \in \bar{p}} A_i$ is either a contractive or a strictly contractive $p-$ semi-cyclic impulsive self-mapping according to Definition 2.2 or Definition 2.4 as stated in the subsequent result:

**Corollary 3.4**. Theorem 3.2 holds, in particular, if $T:\bigcup_{i\in\bar{p}} A_i \to \bigcup_{i\in\bar{p}} A_i$ is a contractive or strictly contractive $p-$ semi-cyclic impulsive self-mapping with $K_i = K \in [0,1)$; $\forall i \in \bar{p}$ being a constant in (2.22) or (2.32) subject to (2.28) and $m:\left(\bigcup_{i\in\bar{p}} A_i\right)\times\left(\bigcup_{i\in\bar{p}} A_i\right) \to \mathbf{R}_{0+}$ being non larger than unity.

Theorem 3.2 also holds if $T:\bigcup_{i\in\bar{p}} A_i \to \bigcup_{i\in\bar{p}} A_i$ is, in particular, a contractive or strictly contractive $p-$ cyclic impulsive self-mapping with $K_i = K \in [0,1)$; $\forall i \in \bar{p}$ being constant and subject to (2.28) and $m:\left(\bigcup_{i\in\bar{p}} A_i\right)\times\left(\bigcup_{i\in\bar{p}} A_i\right) \to \mathbf{R}_{0+}$ being non larger than unity. In this case, the limiting set (3.3) only contains best proximity points, that is, it is of the form $(z_1, z_2, ..., z_p)$. □